\documentclass[11pt]{article}

\usepackage[utf8]{inputenc}
\usepackage[T1]{fontenc}

\usepackage[margin=1in]{geometry}   % set the margins to 1in on all sides
\usepackage{setspace}               % set spacing
\usepackage{lscape}
\usepackage[style=numeric]{biblatex}
\usepackage{graphicx}               % to include figures
\usepackage[dvipsnames]{xcolor}
\usepackage{forest}
\usetikzlibrary{arrows}

\usepackage{amsmath}                % great math stuff
\usepackage{amsfonts}               % for blackboard bold, etc
\usepackage{amsthm}                 % better theorem environments
\usepackage{mathrsfs}
\usepackage{calrsfs}
\usepackage{thmtools}
\usepackage{thm-restate}

\usepackage{hyperref}
\usepackage{cleveref}
\usepackage{caption}
\usepackage{subcaption}
\usepackage{amssymb }
\usepackage{verbatim}
\usepackage{tikz}
\usetikzlibrary{arrows.meta}
\usepackage{comment}
\usepackage{multicol}

\usepackage{import}
\usepackage{standalone}
\usepackage{subfiles}

\newtheorem{thm}{Theorem}[section]

\newtheorem{lem}[thm]{Lemma}

\newtheorem{cor}[thm]{Corollary}
\newtheorem{conj}[thm]{Conjecture}
\newtheorem{defi}[thm]{Definition}

\newcommand{\CC}{\mathcal{C}}		    % for great circle	
\newcommand{\HH}{\mathcal{H}}		    % for hemisphere	
\date{October 2021}

\title{Great-circle Tree Thrackles}

\author{Karen Collins\footnote{Professor of Mathematics, Wesleyan University}, Cleo Roberts\footnote{Independent Researcher (This paper comes from research done while Ms. Roberts was a graduate student at Wesleyan University)}}

\begin{document}

%\pagenumbering{roman}       %Begin numbering with small Roman numerals

\maketitle
\begin{abstract}
A thrackle is a graph drawing in which every pair of edges meets exactly once. Conway's Thrackle Conjecture states that the number of edges of a thrackle cannot exceed the number of its vertices. Cairns et al (2015) prove that the Thrackle Conjecture holds for great-circle thrackles drawn on the sphere. They also posit that Conway's Thrackle Conjecture can be restated to say that a graph can be drawn as a thrackle drawing in the plane if and only if it admits a great-circle thrackle drawing. We demonstrate that the class of great-circle thrackleable graphs excludes some trees. Thus the informal conjecture from Cairns et al (2015) is not equivalent to the Thrackle Conjecture.
\end{abstract}

\pagebreak

\pagebreak

\doublespacing

\pagebreak
\pagenumbering{arabic}

\section{Introduction}

John Conway introduced the word \emph{thrackle} to the world as a Scottish fishing term. As a teenager on vacation in Scotland with his family, Conway encountered a fisherman holding a tangled fishing line. As Conway told it, the fisherman called his line "thrackled." The term stuck with him and later, as a young mathematician, Conway used the word "thrackle" to describe a type of graph drawing in which every pair of edges meets exactly once. He suggested that the number of edges in a thrackle must be bound above by the number of vertices in the thrackle. Since Conway introduced the term to the mathematical community, no one has yet been able to confirm that "thrackle" is a regional term, let alone one having anything to do with fishing. Despite its apocryphal origins, the thrackle and Conway's conjecture live on. In fact, there is still an unclaimed cash prize offered for a proof that Conway's Thrackle Conjecture is true \cite{mitchell_orourke}.

Beginning at the conference where Conway posed the Thrackle Conjecture, Douglas Woodall \cite{woodall} provided the foundation for most of the research on thrackles that followed. Combinatorial attempts to prove the Thrackle Conjecture have relied not only on Woodall's work, but on the research of Lovasz, Pach, and Szegedy \cite{lps}. This research focuses on calculating the upper bound on the number of edges in a thrackleable graph. Fulek and Pach \cite{fulek_pach1}, \cite{fulek_pach2} and Goddyn and Xu \cite{goddyn_xu} continued this work. Xu \cite{xu_2} has published the latest upper bound. 

In doing his initial research about thrackles, Woodall \cite{woodall} demonstrated that the Thrackle Conjecture is true for graphs that can be drawn as thrackles using straight lines. Although he never worked explicitly on thrackles, Paul Erdos \cite{erdos} also provided a design-theory proof that straight-line thrackles obey the Thrackle Conjecture. Woodall also showed that the Thrackle Conjecture does not hold on the torus. On the sphere, however, Cairns, Koussas, and Nikolayevsky \cite{ckn} showed that the Thrackle Conjecture is true for thrackles drawn using arcs of great-circles on the sphere. Since great-circle thrackleable graphs are thrackleable in the plane, Cairns et al's result left open the possibility of proving the Thrackle Conjecture by showing all graphs that are thrackleable in the plane can be drawn as great circle thrackles. However, in this paper we show that there exists a class of thrackleable trees that cannot be drawn as great-circle thrackles. Building on the work of Cairns, Koussas, and Nikolayevsky \cite{ckn}, we classify these trees and begin to characterize great-circle thrackleabe trees.

In the remainder of Section 1 we introduce the concept of thrackles generally, and great-circle thrackles in particular. In Section 2 we present some basic facts about thrackles and known characteristics of great-circle thrackles. Finally, in Section 3 we demonstrate that some trees are not great-circle thrackleable, the main result of the paper.

We begin with a formal definition of a thrackle.

\begin{defi} A \emph{\textbf{thrackle}} is a drawing of a simple, undirected graph in which 
\begin{enumerate}%[)]
\item 
no curve crosses itself and
\item
every pair of curves intersects at exactly one point
\end{enumerate}
If a graph can be drawn as a thrackle, it is called \emph{\textbf{thrackleable}}.
\end{defi}

\begin{figure}[h]
\begin{center}
\begin{tikzpicture}
\filldraw (2,0) circle (2.5pt) node[below] {$v_3$};
\filldraw (1,0) circle (2.5pt) node[below] {$v_2$};
\filldraw (0,0) circle (2.5pt) node[below] {$v_1$};
\filldraw (-1,0) circle (2.5pt) node[below] {$v_0$};
\draw(2,0)--(1,0)--(0,0)--(-1,0);

\draw[-Latex](2.5,0.5)--(3.5,0.5);

\filldraw (5,1) circle (2.5pt) node[left] {$v_3$};
\filldraw (6,0) circle (2.5pt) node[right] {$v_2$};
\filldraw (5,0) circle (2.5pt) node[left] {$v_1$};
\filldraw (6,1) circle (2.5pt) node[right] {$v_0$};
\draw(5,1)--(6,0)--(5,0)--(6,1);
\end{tikzpicture}
\end{center}
\caption{Planar and Thrackle Drawings of 3-Path}
\label{path_3}
\end{figure}
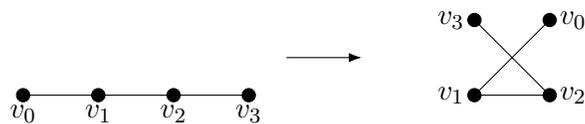

Figures \ref{c3} and \ref{path_3} demonstrate examples of thrackles.

% Figure 1
\begin{figure}[h]
\begin{center}
\begin{tikzpicture}
\filldraw (0,0) circle (2.5pt) node[below] {$v_0$};
\filldraw (-.6,1) circle (2.5pt) node[above left] {$v_1$};
\filldraw (.6,1) circle (2.5pt) node[above right] {$v_2$};
\draw(0,0)--(.6,1)--(-.6,1)--(0,0);
\end{tikzpicture}
\end{center}
\caption{$3$-Cycle}
\label{c3}
\end{figure}

% Figure 2

Conway made the following conjecture about thrackles at a conference in 1969. Despite its age and the many early advances by Douglas Woodall \cite{woodall}, the Thrackle Conjecture remains open.

\begin{conj}[Conway's Thrackle Conjecture, 1969]
Every thrackleable graph contains at least as many vertices as edges.
\end{conj}

\pagebreak

In 2015, Cairns, Koussas, and Nikolayevsky developed a novel approach to the Thrackle Conjecture. They showed that the Conjecture is true for \emph{standard great-circle thrackles}. 

\begin{defi}[Cairns, Koussas, and Nikolayevsky \cite{ckn}]
A \emph{\textbf{great circle thrackle}} is a thrackle drawn on the sphere, whose vertices are represented as points on the sphere and whose edges are represented by the arcs of great circles. A great-circle thrackle is said to be drawn in \emph{\textbf{general position}} if no three vertices lie on the same great circle and no two vertices are antipodal.
\end{defi}

\begin{figure}[h] % Figure 3 
\begin{center} 
\includegraphics[scale=0.5]{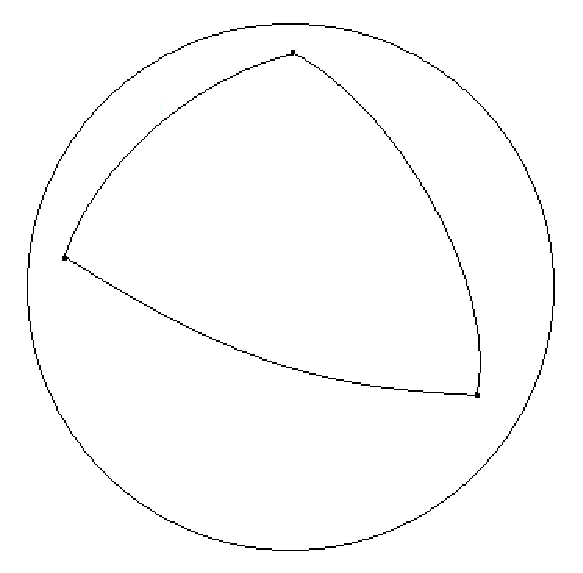} 
\includegraphics[scale=0.5]{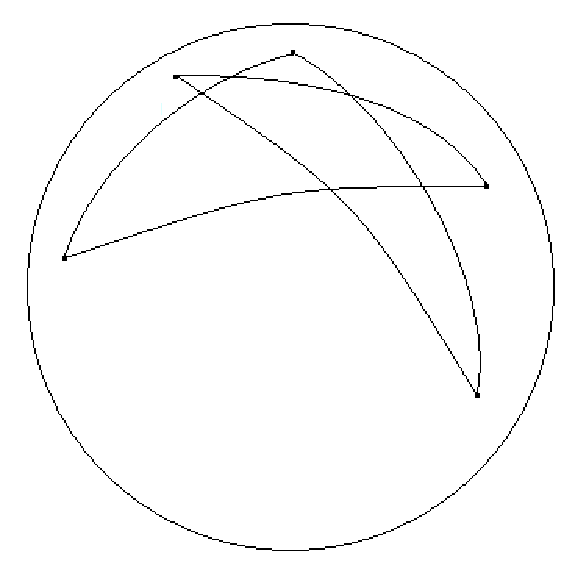} 
\caption{Great-circle Thrackles of the $3$-cycle (left) and $5$-cycle (right)}
\label{great_circle_thrackle}
\end{center} 
\end{figure}

In Figure \ref{great_circle_thrackle}, we show great-circle thrackles of the $3$-cycle and the $5$-cycle, both in general position.

Cairns, Koussas, and Nikolayevsky \cite{ckn} put the further restrictions on great-circle thrackles that they must contain no leaves and must be connected. A great-circle thrackle of a connected graph that contains no leaves and is in general position is a \emph{standard great-circle thrackle}. In addition, they define a the notion of a \emph{crossing orientation} on edges that meet. As we are concerned with trees, we consider connected great-circle thrackles in general position, but we do not require the absence of leaves. Moreover, we do not consider crossing orientation in this analysis.

\begin{thm}[Cairns, Koussas, and Nikolayevsky \cite{ckn}] 
Let $G$ be a standard great circle thrackleable graph. Then $G$ contains at most as many edges as vertices. 
\label{ckn_thm}
\end{thm}

Note that if a great-circle thrackle in general position contradicted the Thrackle Conjecture, it would necessarily contain a standard great-circle thrackle. Therefore, the above theorem holds for any great-circle thrackle in general position. Moreover, Cairns, Koussas, and Nikolayevsky \cite{ckn} demonstrated that every great-circle thrackleable graph is thrackleable in the plane. This result was exciting because it suggested that the Thrackle Conjecture could be reduced to showing that every thrackleable graph can be drawn as a great-circle thrackle. However, we contend that this is not true. We show that there exists a class of trees that are thrackleable in the plane, but cannot be drawn as great circle thrackles. Thus the Thrackle Conjecture must be proved through other means.

\section{Facts about Great-circle Thrackles}

One subclass of great-circle thrackleable graphs has been known for decades. Because great circles on the sphere are analogous to straight lines in the plane, all \emph{straight-line thrackleable} graphs are great-circle thrackleable. In fact, Woodall \cite{woodall} classified straight-line thrackleable graphs, including \emph{caterpillars} shortly after the conference where Conway presented his conjecture.

\begin{defi}
A \emph{\textbf{straight-line thrackle}} is a thrackle whose arcs are drawn as straight line segments.
\end{defi}

\begin{figure}[h]
\begin{center}
\begin{tikzpicture}
% 3-star
\filldraw (-2,-.67) circle (2.5pt) node[below] {$v_0$};
\filldraw (-2,0) circle (2.5pt) node[above left] {$v_1$};
\filldraw (-2.6,-1) circle (2.5pt) node[below left] {$v_2$};
\filldraw (-1.4,-1) circle (2.5pt) node[below right] {$v_3$};
\draw(-2,0)--(-2,-.67);
\draw(-2.6,-1)--(-2,-.67);
\draw(-1.4,-1)--(-2,-.67);

% 5-cycle
\filldraw (2,-1) circle (2.5pt) node[below] {$v_0$};
\filldraw (2.951,-.309) circle (2.5pt) node[right] {$v_{3}$};
\filldraw (2.588,.809) circle (2.5pt) node[above right] {$v_{1}$};
\filldraw (1.412,.809) circle (2.5pt) node[above left] {$v_{4}$};
\filldraw (1.049,-.309) circle (2.5pt) node[left] {$v_{2}$};
\draw(2,-1)--(2.588,.809)--(1.049,-.309)--(2.951,-.309)--(1.412,.809)--(2,-1);
\end{tikzpicture}
\end{center}
\caption{$3$-Star (left) and $5$-Cycle (right)}
\label{straight_line}
\end{figure}
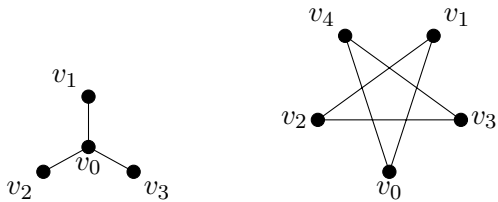

In Figure \ref{straight_line}, we show two examples of straight-line thrackles.

% Figure 4

\begin{figure}[h]
\begin{center}
\begin{tikzpicture}
\filldraw (0,1) circle (2.5pt) node[left] {$v_0$};
\filldraw (-.6,2) circle (2.5pt) node[left] {$v_1$};
\filldraw (.6,2) circle (2.5pt) node[right] {$v_2$};
\filldraw (-.6,0) circle (2.5pt) node[below] {$v_4$};
\filldraw (.6,0) circle (2.5pt) node[below] {$v_3$};
\filldraw (1,1) circle (2.5pt) node[below left] {$v_5$};
\filldraw (1.6,2) circle (2.5pt) node[right] {$v_7$};
\filldraw (1.6,0) circle (2.5pt) node[below] {$v_8$};
\filldraw (2,1) circle (2.5pt) node[below left] {$v_9$};
\filldraw (2.6,2) circle (2.5pt) node[right] {$v_{11}$};
\filldraw (3,1) circle (2.5pt) node[right] {$v_{12}$};
\filldraw (2.6,0) circle (2.5pt) node[below] {$v_{13}$};
\draw(0,1)--(-.6,2);
\draw(0,1)--(.6,2);
\draw(0,1)--(-.6,0);
\draw(0,1)--(.6,0);
\draw(1,1)--(1.6,2);
\draw(1,1)--(1.6,0);
\draw(2,1)--(2.6,2);
\draw(2,1)--(3,1);
\draw(2,1)--(2.6,0);
\draw(0,1)--(1,1)--(2,1);
\end{tikzpicture}
\end{center}
\caption{Caterpillar}
\label{caterpillar}
\end{figure}
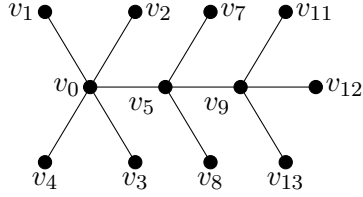

One important straight-line thrackleable graph is the \emph{caterpillar}. Woodall \cite{woodall} showed that any straight-line thrackleable tree must be a caterpillar. We give an example in Figure \ref{caterpillar}.

\begin{defi} A \emph{\textbf{caterpillar}} is a tree in which every vertex is adjacent to at most two interior vertices. 
\end{defi}

% Figure 5

\pagebreak

\begin{thm}[Woodall \cite{woodall}] A finite graph $G$ is straight-line thrackleable if and only if either \begin{enumerate}
        \item $G$ contains an odd cycle and every vertex of $G$ is adjacent to a vertex of the cycle or
        \item $G$ is a disjoint union of caterpillars
    \end{enumerate} \label{thm_slthrackles} \end{thm}
    
To discuss the broader class of great-circle thrackleable trees, we need to define some important terms. We start with \emph{long} and \emph{short} edges.

\begin{defi}[Cairns, Koussas, and Nikolayevsky \cite{ckn}] We represent the great circle containing edge $e$ as $\CC(e)$ and denote the positive and negative hemispheres of $\CC(e)$ by $\HH^{+}_{e}$ and $\HH^{-}_{e}$, respectively. If $e$ is shorter than $\pi$ we call it a \emph{\textbf{short edge}} and if $e$ is longer than $\pi$ we call it a \emph{\textbf{long edge}}. \end{defi}

\begin{figure}[h] % Figure 6 
\begin{center} 
\includegraphics[scale=0.5]{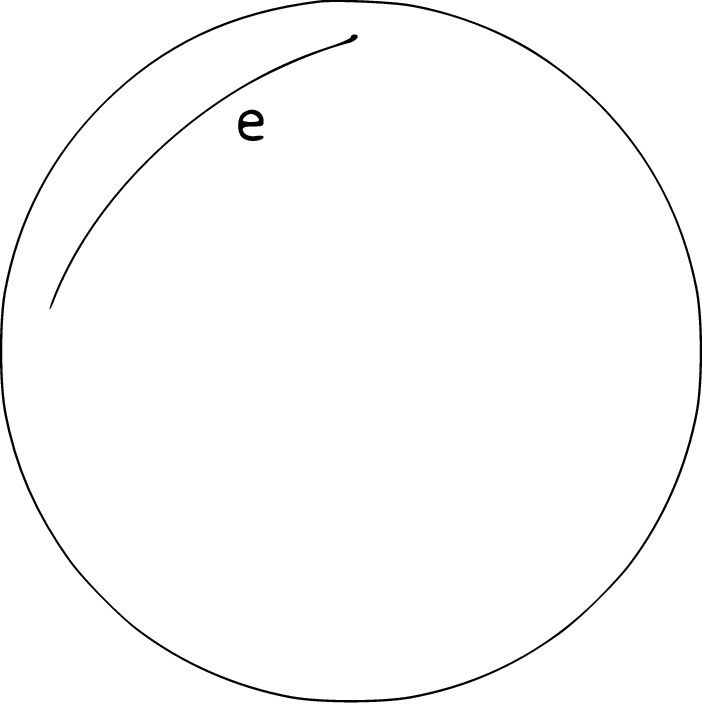} 
\includegraphics[scale=0.5]{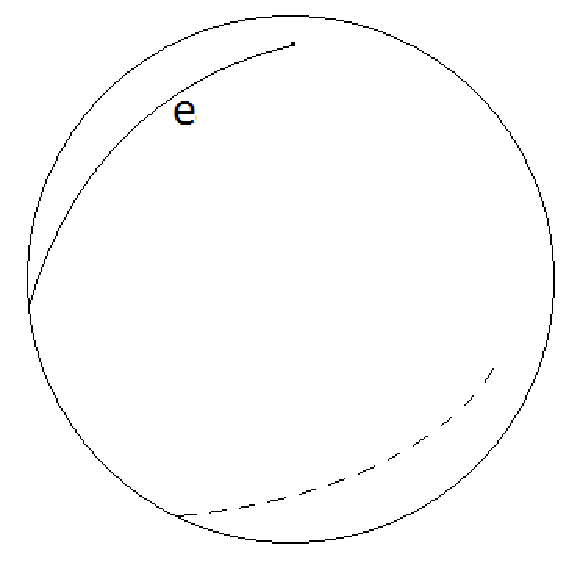} 
\caption{Short (left) and Long (right) Edges}
\label{short_long_edge} 
\end{center} 
\end{figure}

In Figure \ref{short_long_edge}, we have illustrations of a short and a long edge.

One restriction on long edges in great-circle thrackles comes almost directly from the definition.

\begin{lem}[Cairns, Koussas, and Nikolayevsky \cite{ckn}] 
In any great-circle thrackle, no two long edges are adjacent. \label{lem_ckn2015_2a}
\end{lem}

In fact, this lemma is holds for any great-circle thrackle in general position. If, in a great-circle thrackle of graph $G$, any two long edges were adjacent, they would meet at both their shared vertex and the antipodal point of their shared vertex, which would contradict $G$ being drawn as a thrackle.

To further characterize edges in a great-circle thrackle, we describe how an edge \emph{reaches} a vertex and how edges and paths can \emph{separate} at a vertex.

\begin{defi}[Cairns, Koussas, and Nikolayevsky \cite{ckn}] Edge $e$ \emph{\textbf{reaches}} vertex $v$ through hemisphere $\HH$ if $e$ is incident to $v$ and in a small neighborhood of $v$, the interior of $e$ is in $\HH$. Let $v$ be a vertex of $G$ with degree at least $3$. Edge $e$ \emph{\textbf{separates}} at $v$ if there exist edges $f, g \in G$ incident to $v$ such that $f$ reaches $v$ through $\HH^{+}_{e}$ and $g$ reaches $v$ through $\HH^{-}_{e}$. 
\end{defi}

\begin{figure}[h] % Figure 7 
\begin{center} 
\includegraphics[scale=0.5]{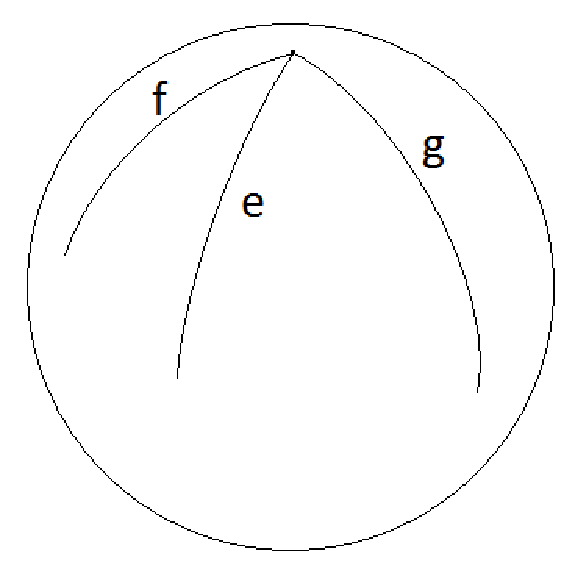} 
\caption{Edge separating at a vertex}
\label{reach_separate} 
\end{center} 
\end{figure}

In Figure \ref{reach_separate}, we see edge $e$ separating at a vertex shared with edges $f$ and $g$.

For the purpose of clarity, we consider only hemispheres bounded by great circles through a specified vertex, $v$. In addition to edges that separate at $v$, we consider paths that separate at $v$.

\begin{defi} A $k$-path $P = e_1\dots e_k$ \emph{\textbf{separates}} at vertex $v$ if either $e_1$ or $e_k$ separates at $v$ and $v$ is either the initial or terminal vertex of $P$.
\end{defi}

\begin{figure}[h] % Figure 8 
\begin{center} 
\includegraphics[scale=0.5]{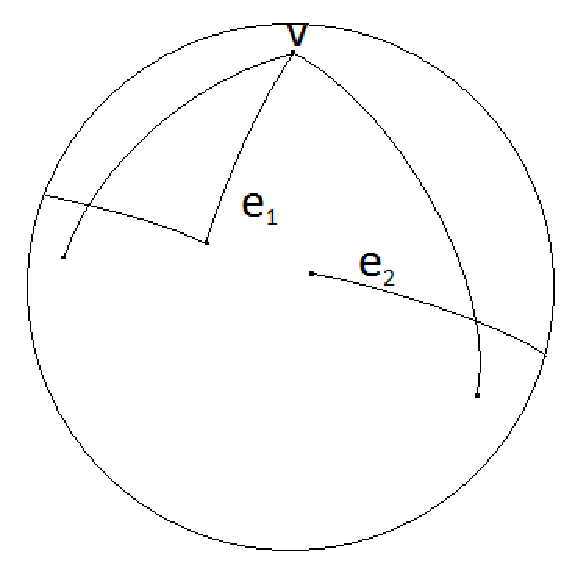} 
\caption{$2$-Path separating at vertex $v$}
\label{path_separate} 
\end{center} 
\end{figure}

In Figure \ref{path_separate}, a $2$-path separates at vertex $v$.

%\pagebreak

\pagebreak

Let us continue to characterize great-circle thrackleable graphs and how adjacent edges relate to the vertices they share. For example, if $v$ is a \emph{pointed vertex} with degree at least three, then at least one edge incident to $v$ separates at $v$.

\begin{defi} If $v$ is a \emph{\textbf{pointed vertex}} in a graph drawing on the sphere, then all of the edges incident to $v$ reach $v$ through a single hemisphere whose boundary is a great circle through $v$. 
\end{defi}

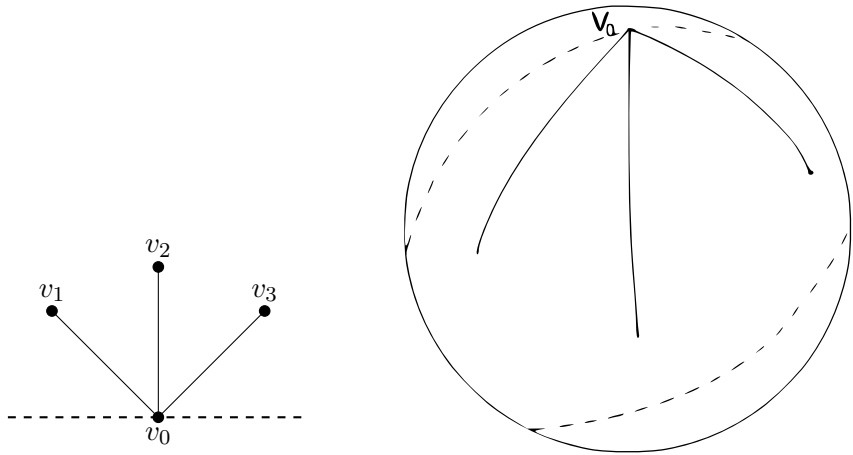
\begin{figure}[h]
\begin{center}
\begin{tikzpicture}[scale=2]
\filldraw (0,0) circle (1pt) node[below] {$v_0$};
\filldraw (-.707,.707) circle (1pt) node[above] {$v_1$};
\filldraw (0,1) circle (1pt) node[above] {$v_2$};
\filldraw (.707,.707) circle (1pt) node[above] {$v_3$};
\draw (0,0)--(-.707,.707);
\draw (0,0)--(0,1);
\draw (0,0)--(.707,.707);
\draw[thick,dashed] (-1,0)--(1,0);
\end{tikzpicture}
\hspace{1cm}
\includegraphics[scale=0.5]{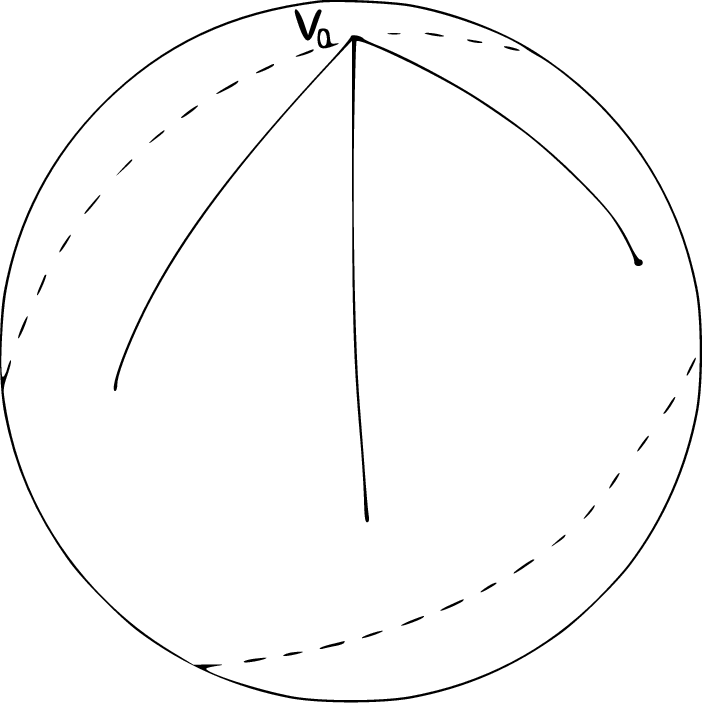} 
\end{center}
\caption{Pointed Vertex}
\label{pointed_vertex}
\end{figure}

Figure \ref{pointed_vertex} provides examples of pointed vertices, $v_0$.

% Figure 9

\begin{lem}Let $G$ be a great-circle thrackleable graph and let $deg(v)=k$ for some $v$ in $G$. Then at least $k-2$ edges separate at $v$. If $v$ is pointed, then exactly $k-2$ edges separate at $v$. \label{lem_robertsseparate} \end{lem} 
\begin{proof} Label the edges incident to $v$ $e_{1}, \dots, e_{k}$. Let $l_i$ be the length of edge $e_i$ and let $U$ be the circular neighborhood of radius $r$ centered at $v$, where $0<r<min\{l_{i}\}_{1 \leq i \leq k}$.

Suppose that $k=3$ and that neither $e_1$ nor $e_2$ separates at $v$. Then $e_2$ and $e_3$ reach $v$ through the same hemisphere of $\CC(e_{1})$, say $\HH_{e_{1}}^{+}$, and $e_1$ and $e_3$ reach $v$ through the same hemisphere of $\CC(e_{2})$, say $\HH_{e_{2}}^{+}$. Thus $e_3$ reaches $v$ through $A = U \cap (\HH_{e_{1}}^{+} \cap \HH_{e_{2}}^{+})$, which is bounded in part by $e_1$ and $e_2$. Because neither $e_1$ nor $e_2$ separates at $v$, the angle between $e_1$ and $e_2$ in $A$ is less than $\pi$. Moreover, $\CC(e_{3})$ passes through the interior of $A$, so $e_1$ and $e_2$ reach $v$ through opposite hemispheres of $\CC(e_{3})$.
    
Now suppose that $k = n$ and that the claim holds if one edge incident to $v$ is removed. Call this removed edge $e_n$. If either all or all but one edge incident to $v$ separates at $v$, then the re-insertion of $e_n$ in $G$ fulfills the claim. We complete the proof by assuming that two edges in $G \setminus e_n$ do not separate at $v$. Denote these edges $e_j$ and $e_l$, where $j,l \in [1, \dots, n-1]$. Re-insert $e_n$ and consider the subgraph $(e_{j} \cap e_{l} \cap e_{n})$ of $G$. By the base case, one of the edges  $e_{j}, e_{l}, e_{n}$ separates at $v$.
    
Suppose that $v$ is pointed. Denote by $\HH_v$ the hemisphere through $v$ through which all edges incident to $v$ reach $v$. Let $\CC(v)$ be the boundary of $\HH_v$. Let $e_1$ and $e_k$ be the edges drawn at the smallest positive angle from $\CC(v)$. Then edges $e_{2}, \dots e_k$ reach $v$ in the portion of $\HH_v$ bounded by $\CC(v)$ and $\CC(e_{1})$. Similarly, edges $e_{1}, \dots e_{k-1}$ reach $v$ in the portion of $\HH_v$ bounded by $\CC(v)$ and $\CC(e_{k})$. Thus exactly $k-2$ edges separate at $v$. \end{proof}

Toward defining the class of great-circle thrackleable trees, we consider edge ajacencies. In Lemmas \ref{lem_gclongseparate} and \ref{lem_robertspath2} we demonstrate some limits on edges and paths that separate at shared vertices.

\begin{lem} Every long edge that separates at a vertex is terminal. \label{lem_gclongseparate} \end{lem} 
\begin{proof}
Let $G$ be a great circle thrackleable graph. Let $v_0 \in G$ be a vertex of degree $k \geq 3$ and $e_1$ an edge that separates at $v$. Let $f$ and $g$ be adjacent to $e$ via $v_0$. Suppose that $e_1$ is drawn long. Toward a contradiction, suppose that $e_2$ is adjacent to $e_1$ via $v_1$. Since $f$ and $g$ are in different hemispheres of $\CC(e_{1})$, $e_2$ must be long to cross them both. Yet this is not possible, by Lemma \ref{lem_ckn2015_2a}. 
\end{proof}

%With Lemma \ref{lem_robertspath2}, we begin to demonstrate limits on tree thrackles.

\begin{lem} Let $G$ be a great circle thrackleable graph and let $v \in G$ be a vertex of degree $k \geq 3$. If all edges incident to $v$ are short, then any path that separates at $v$ contains at most two edges. \label{lem_robertspath2} \end{lem} 
\begin{proof} Let path $P$ separate at degree-$k$ vertex $v_0$, with $e_1 \subseteq P$ incident to $v_0$. Suppose that the length of $P$ is greater than one, with $e_1$ and $e_2$ incident to $v_1$ and $e_2$ incident to $v_2$. Since $e_1$ separates at $v_0$, there exist some $f$ and $g$ that reach $v_0$ through opposite hemispheres of $\CC(e_1)$. By Lemma \ref{lem_gclongseparate}, $e_1$ is drawn short. Because $f$ and $g$ are in opposite hemispheres of $\CC(e_1)$ and $v_1$ lies on $\CC(e_1)$, $e_2$ must be long to cross both $f$ and $g$. Without loss of generality, suppose that $e_2$ reaches $v_1$ through the hemisphere of $\CC(e_{1})$ in which $f$ reaches $v_0$, say $\HH_{e_{1}}(f)$. Then $v_2$ lies in $\HH_{e_{1}}(g)$, between $g$ and $e_1$, as in Figure \ref{lem_robertspath2}.
    
\begin{figure}[h] % Figure 10 
\begin{center} 
\includegraphics[scale=0.5]{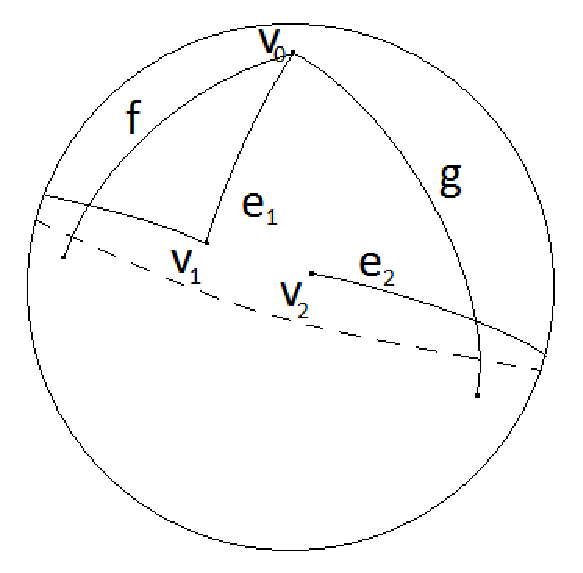} 
\caption{Illustration of path $P$ in Lemma \ref{lem_robertspath2}}
\label{lem_robertspath2} 
\end{center} 
\end{figure}

Now suppose toward a contradiction that $P$ contains a third edge, $e_3$. Because $v_2$ lies in $\HH_{e_{1}}(g)$, but $e_3$ must cross both $f$ and $g$, $e_3$ must be long, which is not possible by Lemma \ref{lem_ckn2015_2a}.
\end{proof}

Theorem \ref{thm_robertsspider2} demonstrates that the class of trees that are great-circle thrackleable properly contains the class of straight-line thrackleable trees. It concerns the \emph{spider} on three legs of length two, which is the minimal example of a tree that is not straight-line thrackleable.

\begin{defi} A \emph{\textbf{spider}} is a tree in which exactly one vertex, $v_0$ has degree at least $3$ and all others have degree at most $2$. The degree of $v_0$ is the number of \emph{\textbf{legs}} of the spider. The \emph{\textbf{length}} of each leg is the length of the path beginning at a degree-$1$ vertex and terminating at $v_0$. 
\end{defi}

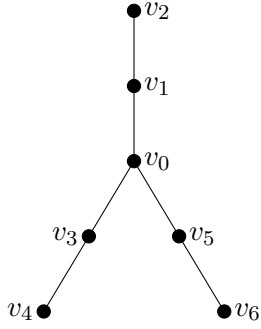
\begin{figure}[h]
\begin{center}
\begin{tikzpicture}
\filldraw (0,1) circle (2.5pt) node[right] {$v_0$};
\filldraw (0,2) circle (2.5pt) node[right] {$v_1$};
\filldraw (0,3) circle (2.5pt) node[right] {$v_2$};
\filldraw (-.6,0) circle (2.5pt) node[left] {$v_3$};
\filldraw (-1.2,-1) circle (2.5pt) node[left] {$v_4$};
\filldraw (.6,0) circle (2.5pt) node[right] {$v_5$};
\filldraw (1.2,-1) circle (2.5pt) node[right] {$v_6$};
\draw(0,1)--(0,2)--(0,3);
\draw(0,1)--(-.6,0)--(-1.2,-1);
\draw(0,1)--(.6,0)--(1.2,-1);
\end{tikzpicture}
\end{center}
\caption{Spider with Three Legs of Length 2}
\label{spider_3legs_length2}
\end{figure}

In Figure \ref{spider_3legs_length2} we show the spider on three legs of length two.

% Figure 11

\pagebreak

\begin{thm} 
The class of great-circle thrackleable graphs properly contains the class of straight-line thrackleable graphs. \label{thm_robertsspider2} 
\end{thm} 
\begin{proof} 
It suffices to show that the spider on three legs of length two can be drawn as a great-circle thrackle. Let us construct $G$, a spider on three legs, each of length two. Let $v$ be its degree-$3$ vertex and let $v$ be pointed. Fix a hemisphere through $v$ containing containing all edges incident to $v$ and call it $S$. Fix a short edge incident to $v$ in $S$ and call it $e_1$. Let $e_1$ be an edge in path $P$. Suppose that $P$ contains a second edge, $e_2$, so that $e_2$ is long and $P$ is drawn in general position. Because $\CC(e_{2})$ is a great circle and since we have the assumption of general position, every point on $\CC(e_{2})$ has distance less than $\pi$ from $v$. Let $max_{e_{2}}$ be the maximum distance from $v$ to any point on $\CC(e_{2})$ within $S$ and let $min_{e_{2}}$ be the minimum distance from $v$ to any point on $\CC(e_{2})$ within $S$. 
        
Draw $f_1$ in $\HH_{e_{1}}^{+} \cap S$ so that it crosses $e_2$ with length $l_{f_{1}}$ greater than $max_{e_{2}}$, but less than $\pi$. Similarly, draw $g_1$ in $\HH_{e_{1}}^{-} \cap S$ so that it crosses $e_2$ with length $l_{g_{1}}$ greater than $max_{e_{2}}$, but less than $\pi$. Let $w_f$ be the endpoint of $f_1$ and $w_g$ the endpoint of $g_1$. So far, this graph is a spider on three legs, two of which have length one. It is also a caterpillar, as we see in Figure \ref{thm_robertsspider2a}.
        
\begin{figure}[h] % Figure 12 
\begin{center} 
\includegraphics[scale=0.5]{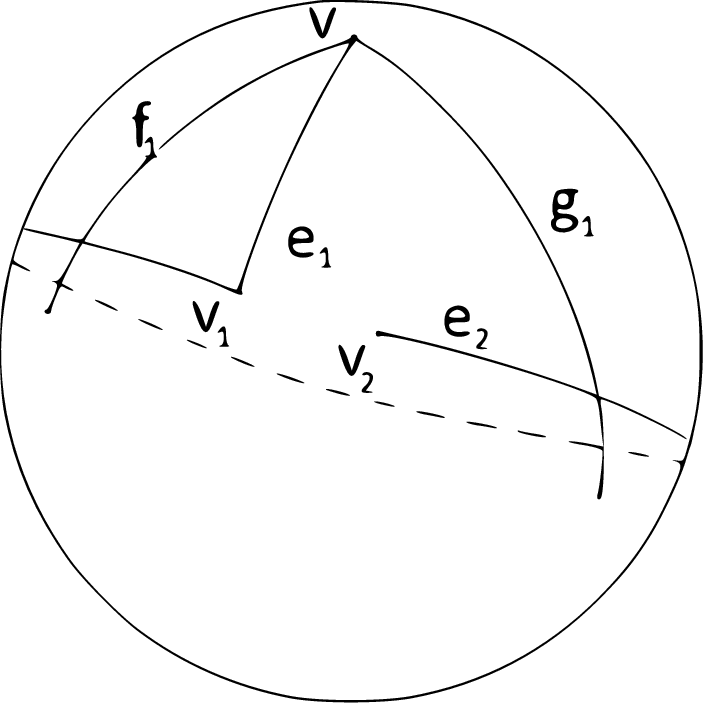} 
\caption{Subgraph of graph $G$}
\label{thm_robertsspider2a} 
\end{center} 
\end{figure}

\pagebreak

Let $p_f$ be a point on $f_1$ with distance from $v$ less than $min_{e_{2}}$ and $p_g$ a point on $g_1$ with distance from $v$ less than $min_{e_{2}}$. Draw $f_2$ as a short arc on the great circle containing $w_f$ and $p_g$ so that $p_g$ is in the interior of $f_2$. Draw $g_2$ similarly, as in Figure \ref{thm_robertsspider2b}.

\begin{figure}[h] % Figure 13 
\begin{center} 
\includegraphics[scale=0.5]{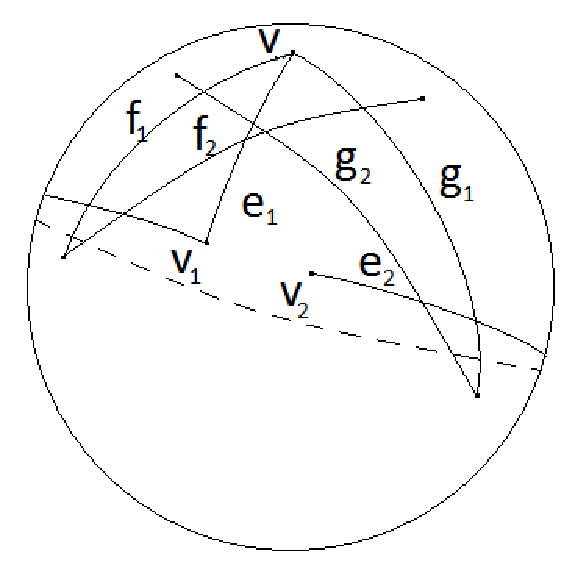} 
\caption{Great-circle Thrackle of Spider on Three Legs of Length Two}
\label{thm_robertsspider2b} 
\end{center} 
\end{figure}

This is a drawing of a spider on three legs, each of length two. \end{proof}

\section{Defining the Class of Great-circle Thrackleable Trees}

Now we are prepared to define a class of trees that contains great-circle thrackleable trees. We call graphs in this class \emph{augmented caterpillars}

\begin{defi}
An \emph{\textbf{augmented caterpillar}}, $G$ is a tree consisting of a \emph{\textbf{spine}}, which is its longest path, and \emph{\textbf{legs}}, which are paths that are edge-disjoint from the spine, but which terminate at internal vertices of the spine. Every vertex in $G$ is at most distance two from an internal vertex of $G$'s spine.
\end{defi}

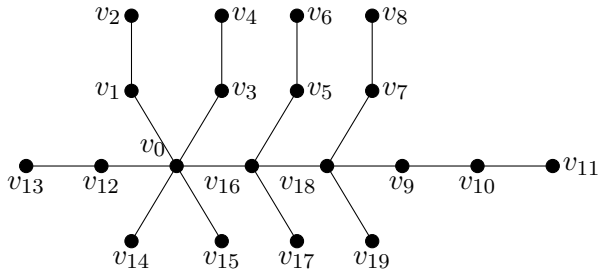
\begin{figure}[h]
\begin{center}
\begin{tikzpicture}
\filldraw (0,1) circle (2.5pt) node[above left] {$v_0$};
\filldraw (-.6,2) circle (2.5pt) node[left] {$v_1$};
\filldraw (-.6,3) circle (2.5pt) node[left] {$v_2$};
\filldraw (.6,2) circle (2.5pt) node[right] {$v_3$};
\filldraw (.6,3) circle (2.5pt) node[right] {$v_4$};
\filldraw (1.6,2) circle (2.5pt) node[right] {$v_5$};
\filldraw (1.6,3) circle (2.5pt) node[right] {$v_6$};
\filldraw (2.6,2) circle (2.5pt) node[right] {$v_{7}$};
\filldraw (2.6,3) circle (2.5pt) node[right] {$v_{8}$};
\filldraw (3,1) circle (2.5pt) node[below] {$v_{9}$};
\filldraw (4,1) circle (2.5pt) node[below] {$v_{10}$};
\filldraw (5,1) circle (2.5pt) node[right] {$v_{11}$};
\filldraw (-1,1) circle (2.5pt) node[below] {$v_{12}$};
\filldraw (-2,1) circle (2.5pt) node[below] {$v_{13}$};
\filldraw (-.6,0) circle (2.5pt) node[below] {$v_{14}$};
\filldraw (.6,0) circle (2.5pt) node[below] {$v_{15}$};
\filldraw (1,1) circle (2.5pt) node[below left] {$v_{16}$};
\filldraw (1.6,0) circle (2.5pt) node[below] {$v_{17}$};
\filldraw (2,1) circle (2.5pt) node[below left] {$v_{18}$};
\filldraw (2.6,0) circle (2.5pt) node[below] {$v_{19}$};
\draw(0,1)--(-.6,2)--(-.6,3);
\draw(0,1)--(.6,2)--(.6,3);
\draw(0,1)--(-.6,0);
\draw(0,1)--(.6,0);
\draw(1,1)--(1.6,2)--(1.6,3);
\draw(1,1)--(1.6,0);
\draw(2,1)--(2.6,2)--(2.6,3);
\draw(2,1)--(3,1)--(4,1)--(5,1);
\draw(2,1)--(2.6,0);
\draw(0,1)--(1,1)--(2,1);
\draw(0,1)--(-1,1)--(-2,1);
\end{tikzpicture}
\end{center}
\caption{Augmented Caterpillar}
\label{augmentedcaterpillar}
\end{figure}

We show an example of an augmented caterpillar in Figure \ref{augmentedcaterpillar}.

% Figure 14

\begin{figure}[h]
\begin{center}
\begin{tikzpicture}
\filldraw (1,0) circle (2.5pt) node[below] {$z$};
\filldraw (2,0) circle (2.5pt) node[below] {$w_1$};
\filldraw (3,0) circle (2.5pt) node[below] {$w_2$};
\filldraw (4,0) circle (2.5pt) node[below] {$w_3$};
\filldraw (0,-.6) circle (2.5pt) node[below right] {$v_1$};
\filldraw (-1,-1.2) circle (2.5pt) node[below right] {$v_2$};
\filldraw (-2,-1.8) circle (2.5pt) node[below right] {$v_3$};
\filldraw (0,.6) circle (2.5pt) node[above right] {$u_1$};
\filldraw (-1,1.2) circle (2.5pt) node[above right] {$u_2$};
\filldraw (-2,1.8) circle (2.5pt) node[above right] {$u_3$};
\draw(1,0)--(2,0)--(3,0)--(4,0);
\draw(1,0)--(0,-.6)--(-1,-1.2)--(-2,-1.8);
\draw(1,0)--(0,.6)--(-1,1.2)--(-2,1.8);
\end{tikzpicture}
\end{center}
\caption{Spider on Three Legs of Length 3}
\label{spider_3legs_length3}
\end{figure}
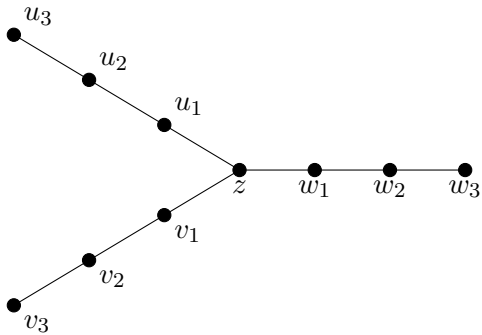

Note that the minimal thrackleable tree which is not an augmented caterpillar is the spider on three legs of length three, shown in Figure \ref{spider_3legs_length3}. 

% Figure 15

If any edge incident to a leaf is removed from this graph, the resulting graph is an augmented caterpillar. If any other edge is removed from this spider, the result is the disjoint union of augmented caterpillars. Thus, by showing that the spider on three legs of length three is not great-circle thrackleable we demonstrate that the class of great-circle thrackleable trees is contained in the class of augmented caterpillars. Consequently, every graph that contains a spider on three legs of length three is cannot be drawn as a great-circle thrackle.

\begin{thm} 
If a tree $G$ is great-circle thrackleable, then $G$ is an augmented caterpillar. 
\label{thm_robertsspider3} 
\end{thm} 
\begin{proof} It suffices to show that a spider on three legs of length three cannot be drawn as a great-circle thrackle. Let $G$ be a spider on three legs of length three. Let $z \in G$ be the degree-$3$ vertex. Label the legs of $G$ as $E = e_{1}e_{2}e_{3}$, $F = f_{1}f_{2}f_{3}$, and $H = h_{1}h_{2}h_{3}$. Let the initial edge in $E$ be $e_1 = \{z, v_1\}$; in $F$, $f_1 = \{z, w_1\}$; and in $H$, $h_1 = \{z, u_1\}$. Let the remaining edges be $e_i = \{v_{i-1}, v_i\}$, $f_i = \{w_{i-1}, w_i\}$, and $h_i = \{u_{i-1}, u_i\}$ for $i = 2,3$, as in Figure \ref{spider_3legs_length3}. 
        
Suppose, toward a contradiction, that $G$ can be drawn as a great circle thrackle. By Lemmas \ref{lem_ckn2015_2a} and \ref{lem_robertspath2}, one of $e_{1}$, $f_{1}$, and $h_1$ is long. Without loss of generality, let $e_1$ be long. By Lemma \ref{lem_gclongseparate}, $e_1$ does not separate at $z$. Again without loss of generality, let $f_1$ separate at $z$, as in Figure \ref{no3legs_step1}. Denote by $l_{e_{1}}$, $l_{f_{1}}$, and $l_{h_{1}}$ be the lengths of $e_{1}$, $f_{1}$, and $h_1$, respectively. To admit a great circle arc containing $w_1$ and interior points of $e_1$ and $h_1$, let $l_{e_{1}} > \pi + max\{l_{f_{1}},l_{h_{1}}\}$. This condition is necessary, but not sufficient to demonstrate that $f_2$ can be drawn to cross both $e_{1}$ and $h_1$, so we further require that $l_{h_{1}} > l_{f_{1}}$. 
        
\begin{figure}[h] % Figure 16 
\begin{center} 
\includegraphics[scale=0.5]{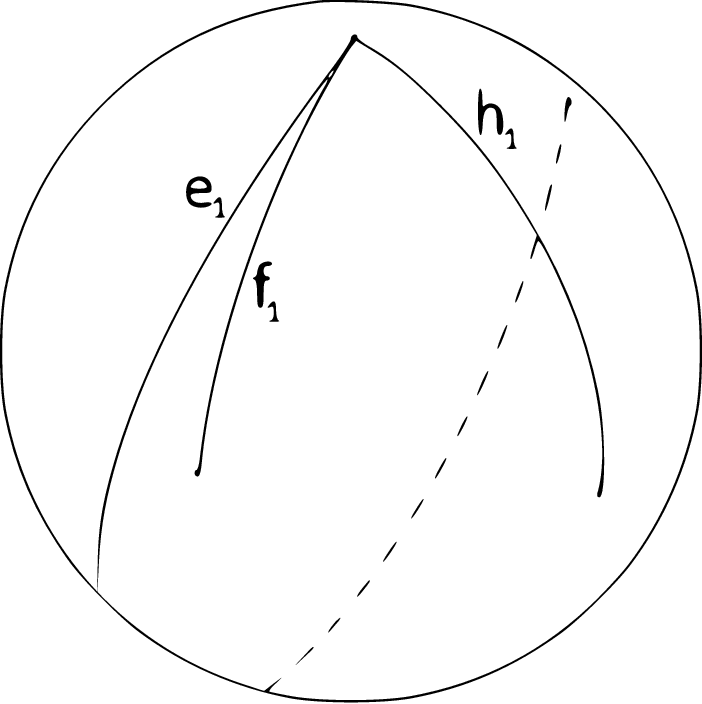}
\includegraphics[scale=0.5]{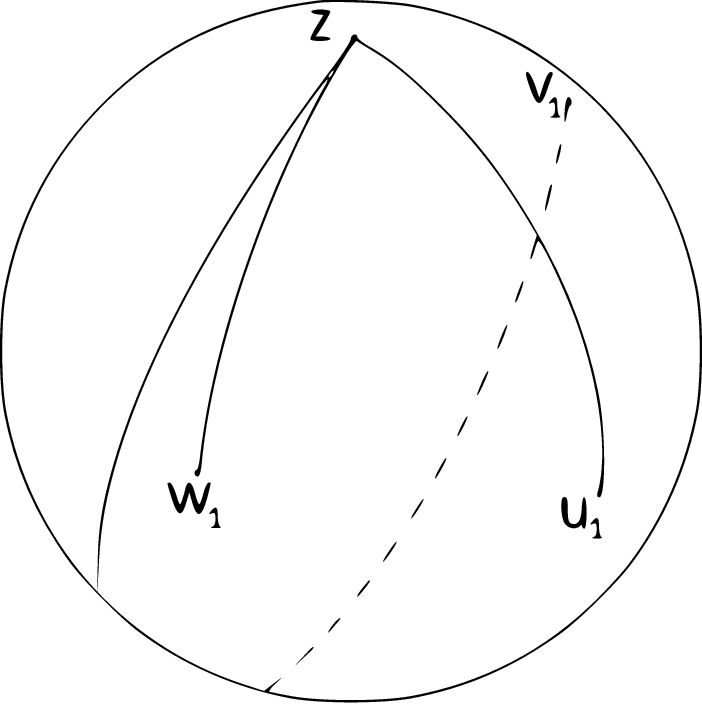} 
\caption{Spider Subgraph of $G$ Containing Three Legs of Length One: Edge labels (left), Vertex labels (right)} \label{no3legs_step1} 
\end{center} 
\end{figure}
        
Although $e_1$ and $h_1$ reach $z$ through different hemispheres of $\CC(f_{1})$, the segment of $e_1$ bounded by $v_1$ and the antipodal point of $z$ is in the same hemisphere of $\CC(f_{1})$ as $h_1$. Draw a short great circle arc $f_2$ from $w_1$ through $h_1$ and $e_1$, terminating at $w_2$. Let $p$ be the point where $f_2$ crosses $h_1$ and $q$ the point where $f_2$ crosses $e_1$, as in Figure \ref{no3legs_step2}.

\begin{figure}[h] % Figure 17 
\begin{center} 
\includegraphics[scale=0.5]{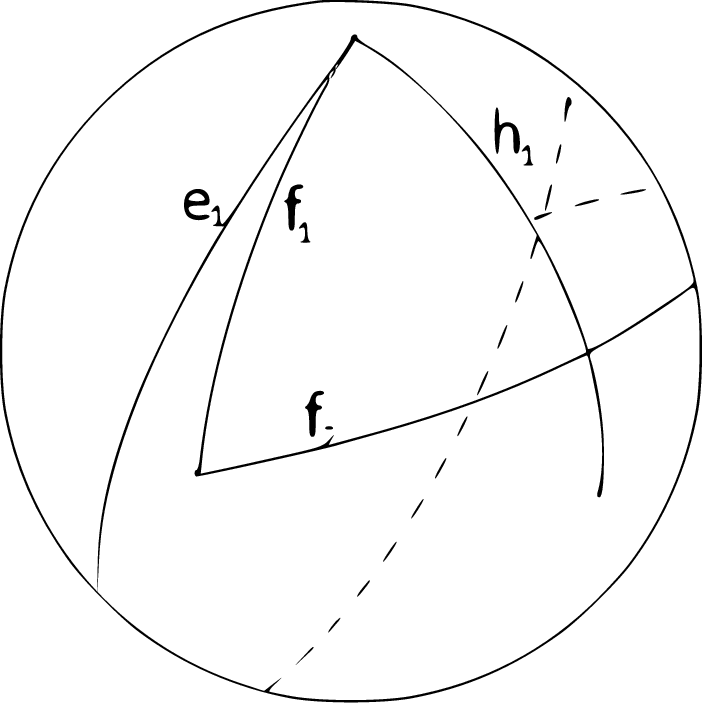}
\includegraphics[scale=0.5]{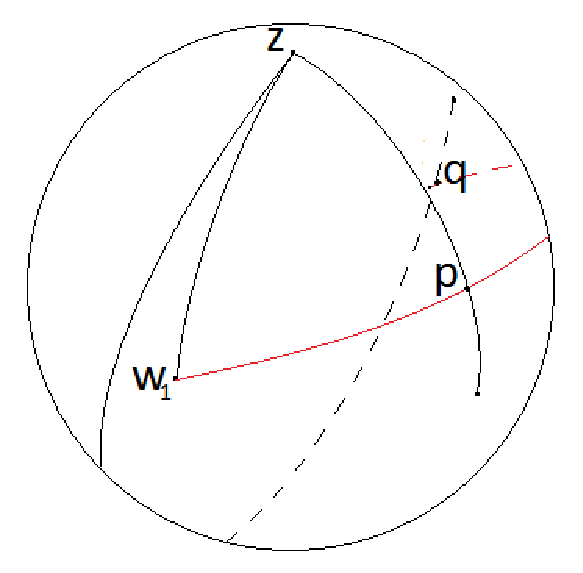} 
\caption{Spider Subgraph of $G$ Containing Two Legs of Length One, One Leg of Length Two: Edge labels (left), Vertex labels (right)} 
\label{no3legs_step2} 
\end{center} 
\end{figure}

\pagebreak
        
Note that the arcs $[z,w_{1}]$, $[w_{1},p]$, and $[p,z]$ form a $3$-cycle. Any edge that is not adjacent to $f_1$, $f_2$, or $h_1$ must cross all three of them. Note that any great circle arc crossing $f_1$, $f_2$, and $g_1$, but not containing more than one interior point of any of these edges, must contain $p$, as in Figure  \ref{no3legs_step3redux} for example.

\begin{figure}[h] % Figure 18 
\begin{center} 
\includegraphics[scale=0.5]{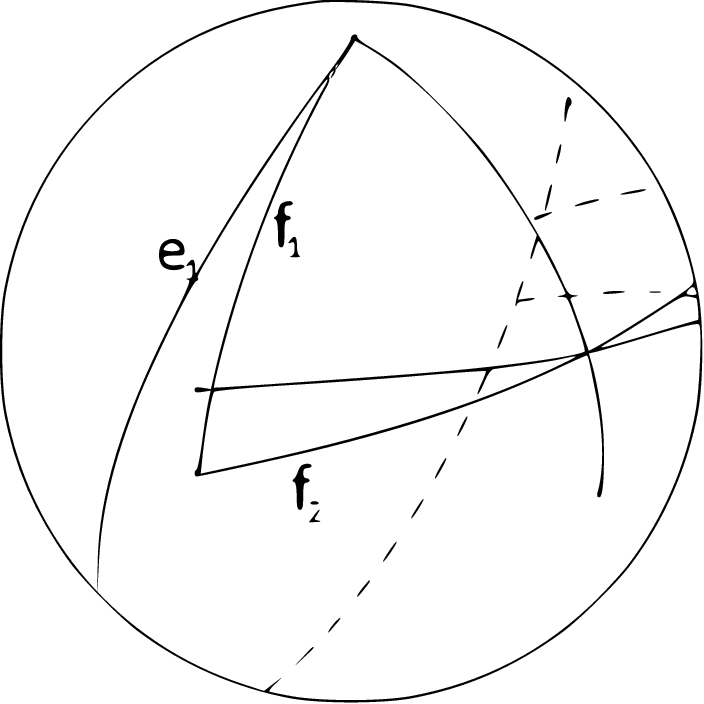}
\includegraphics[scale=0.5]{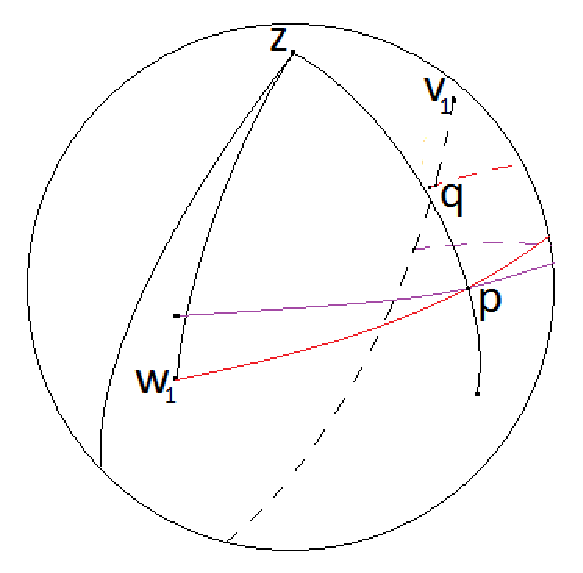} 
\caption{Failure of the Spider on Three Legs of Length Three: Edge labels (left), Vertex labels (right)} \label{no3legs_step3redux} 
\end{center} 
\end{figure}

Recall that $v_1$ lies outside the $3$-cycle formed by $[z,w_{1}]$, $[w_{1},p]$, and $[p,z]$. Suppose that $e_2$ contains both $p$ and an interior point of $f_1$. Then $e_2$ crosses $f_1$ in one hemisphere of $f_2$ and meets $e_1$ in the other hemisphere. In other words, $q$ lies between $v_1$ and the intersection of $e_1$ and $e_2$. This situation is impossible, because the intersection of $e_1$ and $e_2$ is $v_1$.
        
Hence the spider on three legs of length three is not great circle thrackleable. The remainder of the proof is immediate from Lemma \ref{lem_robertspath2} and Theorem \ref{thm_robertsspider2}
\end{proof}

This theorem demonstrates that great-circle thrackleable trees are not the key to proving the Thrackle Conjecture. Moreover, it gives us the following corollary.

\begin{cor}
No graph containing the spider on three legs of length three can be drawn as a great-circle thrackle.
\end{cor}

To conclude, we present the following two conjectures, which may further specify the class of trees that are great-circle thrackleable. The authors suspect that reintroducing the notion of crossing orientation to the discussion of trees may be useful in either proving them or finding counterexamples. The challenge comes from making the notion of crossing orientation well-defined in the absence of closed loops.

\begin{conj} The class of great circle thrackleable trees is properly contained in the class of augmented caterpillars. \end{conj}

\begin{conj} Let $G$ be an augmented caterpillar drawn as a a great circle thrackle and let $v$ and $v'$ be adjacent in the spine of $G$. If some $2$-path separates at $v$, then \emph{no} $2$-path separates at $v'$. \end{conj}

\section{Conclusion}

Despite the many efforts to prove the Thrackle Conjecture (or find a counterexample) over the years, it remains an open problem. The work of Cairns, Koussas, and Nikolayevsky \cite{ckn} showed that the Thrackle Conjecture is true for great-circle thrackleable graphs. Moreover, their work suggested a potential way to prove the Thrackle Conjecture once and for all: show that all graphs that are thrackleable in the plane can be drawn as great circle thrackles. However, this endeavor fails because the class of augmented caterpillars cannot be drawn as great-circle thrackles.

\pagebreak
        
\printbibliography

\end{document}